\documentclass[11pt,reqno]{amsart}
\usepackage[T1]{fontenc}
\usepackage{lmodern}
\usepackage{microtype}
\usepackage[margin=1.1in]{geometry}
\usepackage{amssymb,mathtools}
\usepackage[colorlinks=true,linkcolor=blue,citecolor=blue,urlcolor=blue]{hyperref}

\newtheorem{theorem}{Theorem}

\theoremstyle{remark}
\newtheorem{remark}[theorem]{Remark}
\newcommand{\SU}{\mathrm{SU}}
\newcommand{\SO}{\mathrm{SO}}
\newcommand{\R}{\mathbb R}
\newcommand{\C}{\mathbb C}
\newcommand{\reg}{\mathrm{reg}}

\DeclareMathOperator{\Int}{int}
\DeclareMathOperator{\Gr}{Gr}
\numberwithin{equation}{section}

\title[Self-intersecting sections of polar actions]{Self-intersecting sections of polar actions\texorpdfstring{\\}{} on simply connected manifolds}
\author{Stephan Wiesendorf}
\date{}
\hypersetup{pdftitle={Self-intersecting sections of polar actions on simply connected manifolds},pdfauthor={Stephan Wiesendorf}}

\begin{document}

\begin{abstract}
We construct an effective polar action of $\SO(3)$ on a closed simply
connected $6$-manifold with a compact $3$-dimensional section that is not
injectively immersed. This gives a negative answer to a question of Grove
and Ziller. The principal isotropy is trivial, and the section is injective
over the principal stratum. Its self-intersection set consists of two
disjoint $2$-spheres in a single singular stratum, where exactly two branches
meet with different tangent spaces. The example is obtained by
collapsing three circle subgroups over the boundary components of a
three-dimensional chamber. The self-intersections arise from a quarter-turn
that lies in one of these circles and also in the group generated by the
wall involutions.
\end{abstract}

\maketitle

\section{Introduction}

An isometric action of a compact Lie group $G$ on a complete Riemannian
manifold $M$ is \emph{polar} if there is an isometric immersion
$\iota\colon\Sigma\to M$ of a complete connected manifold which meets
every orbit orthogonally. Such an immersion is called a \emph{section}.
We use the convention of Grove and Ziller \cite{GZ12}, allowing immersed
sections but excluding a multiple covering of another section. Thus the
question of injectivity concerns actual self-intersections, rather than
the choice of a parametrization. Some authors reserve the word polar for
actions with embedded sections and use \emph{locally polar} for this
broader definition; see \cite[p.~3, footnote~2]{Gor25}.

Immediately after Definition~1.1 in \cite{GZ12}, Grove and Ziller ask
whether section immersions must be injective, at least when the ambient
manifold is simply connected. Gorodski and Kollross repeat the
self-intersection question in the introduction to \cite{GK16}, attributing
it to Ziller. It also appears as an open problem in the 2025 revision of
Gorodski's lecture notes \cite[p.~3, footnote~2]{Gor25}. To the best of our knowledge, the question has remained open until now.

The answer is positive in several special cases. For group actions
of cohomogeneity one, sections are injectively immersed
\cite[Theorem~6.1]{AA93}. Sections of asystatic actions, i.e. actions whose isotropy representation on the tangent space to a principal orbit has no nonzero fixed vectors, are properly
embedded \cite{GK16}.

For general singular Riemannian foliations with sections, however,
the author \cite{Wie15} constructed examples on spheres with
one-dimensional sections that can have self-intersections and need
not be mutually isometric. Non-isometric sections cannot occur for
polar group actions, since $G$ acts transitively on the set of
sections. The example below shows, however, that sections of polar
group actions need not be injectively immersed, even when the
ambient manifold is closed and simply connected.

\begin{theorem}\label{thm:main}
There is a closed simply connected Riemannian $6$-manifold $M$ with an
effective polar action of $\SO(3)$ of cohomogeneity three and a compact
connected section $\iota\colon\Sigma^3\to M$ such that:
\begin{enumerate}
\item the principal isotropy is trivial and every isotropy group is
connected;
\item $\iota$ is injective on $\iota^{-1}(M_{\reg})$, where \(M_{\reg}\) is the principal stratum;
\item some points of $M$ have exactly two preimages under $\iota$, and
the two tangent images at each such point are different.
\end{enumerate}
In particular, $\iota$ has no nontrivial subcover section and is not
injective.
\end{theorem}

The construction is a special case of the chamber reconstruction theorem
of Grove and Ziller \cite[Theorem~3.3]{GZ12}, which also provides a smooth
invariant polar metric. We describe such a metric explicitly to explain the geometry near the singular orbits and the
self-intersections of the section. In our example, the group generated by the wall
reflections contains a quarter-turn in one of the wall isotropy circles.
At this wall, only the corresponding half-turn is used to glue the
section. The quarter-turn therefore maps points of the section to
distinct points with the same image in $M$.

The original construction used $\SU(3)$, which made all orbits simply
connected. As suggested by the referee, the same circle subgroups and
wall involutions already lie in $\SO(3)$. Replacing the acting group
by $\SO(3)$ lowers the dimension of the ambient manifold from eleven
to six. The principal orbits are no longer simply connected, but the
ambient manifold still is, as we verify below.

\medskip
\noindent\textbf{Acknowledgements.}
I thank the referee for a careful reading and for suggesting the use
of $\SO(3)$, which substantially reduces the dimension of the example.

ChatGPT, an AI
assistant developed by OpenAI, was used during the development of the example,
for literature checks, and for assistance with drafting. The author assumes full
responsibility for the entire content of the article.

\section{The ambient manifold and its metric}

\subsection{Three circles}
We begin with the subgroups that will determine the singular orbits.
Let $G=\SO(3)$, and denote by $R_a(\theta)$ the rotation through angle
$\theta$ about an oriented unit vector $a\in\R^3$. We choose
\[
 a_1=e_1,\qquad a_2=\frac{e_1+e_2}{\sqrt2},\qquad a_3=e_3,
\]
and let
\[
 K_j=\{R_{a_j}(\theta)\mid\theta\in\R\},\qquad j=1,2,3.
\]
We identify each of these circles with $S^1 \subset \mathbb C$ by
the character $\chi_j\colon K_j\to S^1$ given by
$\chi_j(R_{a_j}(\theta))=e^{i\theta}$. In particular, the half-turn
$s_j=R_{a_j}(\pi)$ is the unique element of order two in $K_j$.
In matrix form these elements are
\[
 s_1=\begin{pmatrix}1&0&0\\0&-1&0\\0&0&-1\end{pmatrix},\quad
 s_2=\begin{pmatrix}0&1&0\\1&0&0\\0&0&-1\end{pmatrix},\quad
 s_3=\begin{pmatrix}-1&0&0\\0&-1&0\\0&0&1\end{pmatrix}.
\]
Their product
\begin{equation}\label{eq:r}
 r=s_1s_2=
 \begin{pmatrix}0&1&0\\-1&0&0\\0&0&1\end{pmatrix}
 =R_{e_3}(-\pi/2)
\end{equation}
satisfies $r^2=s_3$ and $s_1rs_1=r^{-1}$. Thus
$W=\langle s_1,s_2,s_3\rangle$ is the dihedral group of order eight,
with elements
\[
 W=\{e,r,r^2,r^3,s_1,rs_1,r^2s_1,r^3s_1\}.
\]
Since an element of $\SO(3)$ belongs to $K_j$ precisely when it fixes
$a_j$, the above list of elements of $W$ gives
\begin{equation}\label{eq:intersections}
 W\cap K_1=\langle s_1\rangle,\qquad
 W\cap K_2=\langle s_2\rangle,\qquad
 W\cap K_3=\langle r\rangle.
\end{equation}
The third intersection is the one we shall use. It contains the
quarter-turn $r$, whereas the subgroup generated by $s_3$ consists
only of $e$ and $r^2$.

\subsection{The chamber and the quotient}
Choose three pairwise disjoint smoothly embedded closed balls
$B_1,B_2,B_3$ in $S^3$, and remove their interiors. The resulting
manifold
\[
 C=S^3\setminus\bigcup_{j=1}^3\Int B_j
\]
has boundary components $F_j=\partial B_j$, each diffeomorphic to
$S^2$. It is connected and simply connected. Indeed, attaching the
three balls again does not change the fundamental group, by the
van Kampen theorem, and recovers $S^3$.

For each boundary component choose a collar embedding
\[
 \kappa_j\colon F_j\times[0,\varepsilon)\longrightarrow C
 \ \ \text{ with } \ \ \kappa_j(y,0)=y.
\]
After decreasing $\varepsilon$, we may assume that the images
$U_j$ of these embeddings are pairwise disjoint. The parameter $t$
in $\kappa_j(y,t)$ runs from the boundary into the complement of the
removed balls. Such collars can be obtained, for example, by
following the inward normal field of an auxiliary metric on $C$.
We shall use $(y,t)$ as coordinates on $U_j$, with the embedding
$\kappa_j$ understood whenever we regard this pair as a point of $C$.

Over the interior of $C$ we take the product with $G$. Over $F_j$
we replace each copy of $G$ by $G/K_j$. More precisely, let
\begin{equation}\label{eq:M}
 M=(G\times C)/\!\sim
\end{equation}
be the quotient in which points with different $C$-coordinates
are never identified and denote by $q\colon G\times C\to M$ the quotient map. For a fixed $c\in\Int C$, the equivalence
class of $(g,c)$ is a singleton. For $c\in F_j$, it is
\[
 \{(gk,c)\mid k\in K_j\}.
\]
Equivalently, $q(g,c)=q(g',c)$ at such a boundary point if and only
if $g^{-1}g'\in K_j$. Since the boundary components are disjoint, each boundary point belongs to exactly one \(F_j\), which determines the subgroup \(K_j\) used in the identification.
Left multiplication in the first factor preserves this relation and induces an action of \(G\) on \(M\), whose orbit space is naturally identified with \(C\) via \(G\cdot q(g,c)\mapsto c\).

To give \(M\) a smooth structure without boundary, we identify the quotient over each collar with a disk bundle over \(F_j\times G/K_j\). In each fiber, the collar parameter becomes the radial coordinate, while the \(K_j\)-direction supplies the angular coordinate. Collapsing the circle at the boundary then corresponds to adding the center of the disk. Put
\(
 D^2_\varepsilon=\{z\in\C\mid |z|<\varepsilon\}\) and
 \(E_j=G\times_{K_j}D^2_\varepsilon.
\)
Here $K_j$ acts on $G\times D^2_\varepsilon$ on the right by
\[
 (g,z)\cdot k=(gk,\chi_j(k)^{-1}z).
\]
We write $[g,z]_j$ for the corresponding class, so that
\begin{equation}\label{eq:disk}
 [gk,z]_j=[g,\chi_j(k)z]_j.
\end{equation}
The projection $[g,z]_j\mapsto gK_j$ makes $E_j$ a smooth open
disk bundle over $G/K_j$. Its zero section consists of the points
$[g,0]_j$.

The neighborhood $N_j=q(G\times U_j)$ of the singular orbits over
$F_j$ is identified with $F_j\times E_j$ by
\begin{equation}\label{eq:collar-model}
 \Phi_j\colon N_j\longrightarrow F_j\times E_j,
 \qquad
 \Phi_j\bigl(q(g,\kappa_j(y,t))\bigr)=(y,[g,t]_j).
\end{equation}
In this formula the nonnegative real number $t$ is viewed as a
complex number. For $t>0$, every point in the punctured disk bundle
has exactly one representation of this form. In fact, if $z\ne0$
and $k=\chi_j^{-1}(z/|z|)$, then
\[
 [g,z]_j=[gk,|z|]_j.
\]
The inverse of \eqref{eq:collar-model} on this set is therefore
\[
 (y,[g,z]_j)\longmapsto
 q\bigl(g\chi_j^{-1}(z/|z|),\kappa_j(y,|z|)\bigr).
\]
At $z=0$ the inverse sends $(y,[g,0]_j)$ to $q(g,y)$, which depends
only on the coset $gK_j$.

To see that \(\Phi_j\) is a homeomorphism also along the zero section, choose an open set \(V\subset G/K_j\) over which the bundle \(G\to G/K_j\) admits a smooth section \(\sigma\colon V\to G\). Every \(g\in G\) with \(gK_j\in V\) can then be written uniquely as \(g=\sigma(b)k\), where \(b=gK_j\) and \(k\in K_j\).

Using this trivialization and the collar coordinates, a tuple \((y,b,k,t)\in F_j\times V\times K_j\times[0,\varepsilon)\) represents the point \((\sigma(b)k,\kappa_j(y,t))\in G\times C\). Under the quotient map \(q\), all such points with fixed \((y,b)\) and \(t=0\) are identified, while no identifications occur for \(t>0\). In these coordinates, \(\Phi_j\) is induced by
\[
(y,b,k,t)\longmapsto (y,b,t\chi_j(k)).
\]
Since \(\chi_j\) identifies \(K_j\) with the unit circle, the last coordinate is the usual polar-coordinate map onto \(D^2_\varepsilon\). After collapsing the circle at \(t=0\), it becomes a homeomorphism onto the disk. Thus the quotient is locally identified with \(F_j\times V\times D^2_\varepsilon\), including the points with \(t=0\).

We give $N_j$ the smooth structure of $F_j\times E_j$ and retain
the product smooth structure on $q(G\times\Int C)$. On their
overlap, the maps above are diffeomorphisms, since $|z|$ and
$z/|z|$ are smooth for $z\ne0$. Different collar neighborhoods
do not overlap. These charts therefore define a smooth structure
on $M$, and the $G$-action is smooth in these charts. The
equivalence relation in \eqref{eq:M} is closed in
$(G\times C)^2$, so the quotient is compact and Hausdorff.
Consequently $M$ is a closed smooth manifold of dimension six.

\subsection{Orbit types and simple connectivity}
The action is free on
\[
 M_{\reg}=q(G\times\Int C)\cong G\times\Int C.
\]
For $c\in F_j$, the isotropy group at $q(g,c)$ is $gK_jg^{-1}$.
Thus the singular set has three connected components
\[
 S_j=q(G\times F_j)\cong F_j\times G/K_j,
\]
each of codimension two. All three circles $K_j$ are conjugate in
$G$. The principal orbits are copies of $\SO(3)$ and the singular
orbits are copies of $\SO(3)/K_j\cong S^2$. In particular, every
isotropy group is connected, and the action is effective because
it is free on the principal stratum.

To compute the fundamental group of $M$, first note that a loop
can, by general position, be perturbed to avoid the three codimension-two submanifolds
$S_j$. The inclusion of $M_{\reg}$ into $M$ therefore induces a
surjection on fundamental groups. The interior of $C$ has the same
homotopy type as $C$, as one sees by pushing the boundary slightly
inward along the collars. Since $C$ is simply connected, we obtain
\[
 \pi_1(M_{\reg})\cong\pi_1(\SO(3))\cong\mathbb Z_2.
\]
For $c_0\in\Int C$, the nontrivial element is represented by
\[
 \theta\longmapsto q(R_{a_1}(\theta),c_0),
 \qquad 0\leq\theta\leq2\pi.
\]
Choose a path $\gamma\colon[0,1]\to C$ from $c_0$ to a point of
$F_1$. The homotopy
\[
 (\theta,s)\longmapsto q(R_{a_1}(\theta),\gamma(s))
\]
contracts this loop at $s=1$, since the whole circle $K_1$ is
collapsed over $F_1$. Thus the surjection from
$\pi_1(M_{\reg})$ has trivial image, and $M$ is simply connected.
Alternatively, one can apply the van Kampen theorem again. The neighborhood $N_j$
retracts onto $F_j\times G/K_j\cong S^2\times S^2$, while
$N_j\cap M_{\reg}\cong G\times F_j\times(0,\varepsilon)$
is connected and induces an isomorphism on fundamental groups
with $M_{\reg}$. Adding $N_1$ kills its fundamental group, and
adding $N_2$ and $N_3$ does not introduce any new generators.

\subsection{A smooth invariant metric}
Choose a metric $h$ on $C$ for which each collar is a product,
\[
 \kappa_j^*h=h_j+dt^2,
\]
where $h_j$ is a metric on $F_j$. Let $Q$ be a bi-invariant metric
on $G$. The elements $X_j\in\mathfrak{so}(3)$ determined by
\[
 \exp(\theta X_j)=R_{a_j}(\theta)
\]
are conjugate, so a single rescaling of $Q$ makes all of them unit
vectors. We use the same notation for their left-invariant vector
fields on $G$, and write $\theta_j=Q(X_j,\cdot)$ and
$Q_j=Q-\theta_j\otimes\theta_j$. Thus $Q_j$ agrees with $Q$ on
the orthogonal complement of $X_j$ and vanishes in the $X_j$
direction.

We start with the product metric $Q+h$ on $G\times\Int C$.
On the $j$th collar we change only the length of the circle
direction and use
\begin{equation}\label{eq:metric}
 Q_j+f_j(t)^2\theta_j\otimes\theta_j+h_j+dt^2.
\end{equation}
Here $f_j$ is smooth, positive for $t>0$, equal to $t$ near zero
and equal to $1$ near the outer end of the collar. Such a function
can be obtained by a smooth interpolation between $t$ and $1$.
Near the outer end, \eqref{eq:metric} agrees with $Q+h$, so these
formulas give a smooth $G$-invariant metric on $M_{\reg}$.

Since $f_j(t)=t$ near zero, the metric on each disk fiber
is Euclidean near its center, with polar-coordinate expression
$dt^2+t^2d\varphi^2$, where \(\varphi\) denotes the angular coordinate. To check that the metric extends smoothly across \(S_j\), observe that \(\theta_j\) is invariant under right multiplication by \(K_j\).
Indeed, $Q$ is bi-invariant and $\operatorname{Ad}(k)X_j=X_j$ for
$k\in K_j$. Since $\theta_j(X_j)=1$, it is a connection form for
$G\to G/K_j$. The tensor $Q_j$ is the pullback of a metric
$\overline Q_j$ on $G/K_j$.

Choose a local section $\sigma\colon V\to G$ as above and put
$\alpha_j=\sigma^*\theta_j$. In the coordinates
$g=\sigma(b)R_{a_j}(\varphi)$ we have
\[
 \theta_j=d\varphi+\alpha_j,\qquad z=t e^{i\varphi}.
\]
The terms involving $t$ and the circle direction can then be
written as
\[
 dt^2+t^2(d\varphi+\alpha_j)^2=|dz+iz\alpha_j|^2.
\]
Thus, near the zero section of the disk bundle, the metric is
\[
 h_j+\overline Q_j+|dz+iz\alpha_j|^2.
\]
This expression is smooth and positive definite at $z=0$, where
it reduces to $h_j+\overline Q_j+|dz|^2$. It gives the required
extension to $M$. The extended metric is \(G\)-invariant and, since \(M\) is compact, complete. On the principal stratum, the $C$-directions
remain orthogonal to the orbits and carry the metric $h$.

\section{The section and its self-intersections}

We construct the section by gluing eight copies of $C$, one for
each element of $W$. In the copy indexed by $w$, we glue the
boundary component $F_j$ to the same boundary component in the
copy indexed by $ws_j$, using the identity map on $F_j$. In other
words, define
\begin{equation}\label{eq:Sigma}
 \Sigma=(W\times C)/\!\approx
 \ \ \text{ with } \ \ (w,c)\approx(ws_j,c)\quad\text{for }c\in F_j.
\end{equation}
We write \([w,c]\) for the class of \((w,c)\) in \(\Sigma\). There are no identifications over $\Int C$.
At a boundary point, the equivalence class consists of exactly
$(w,c)$ and $(ws_j,c)$, because the boundary components are disjoint
and $s_j^2=e$.

The smooth structure near a glued boundary is described by a signed
collar coordinate. For $y\in F_j$ and $-\varepsilon<t<\varepsilon$,
put
\begin{equation}\label{eq:seam}
 \eta_{w,j}(y,t)=
 \begin{cases}
 [w,\kappa_j(y,t)]&\text{if }t\geq0,\\
 [ws_j,\kappa_j(y,-t)]&\text{if }t\leq0.
 \end{cases}
\end{equation}
At \(t=0\), the points \([w,y]\) and \([ws_j,y]\) are identified by the defining relation for \(\Sigma\), hence the two formulas agree at $t=0$. These product neighborhoods,
together with the interiors of the copies of $C$, make $\Sigma$
a closed smooth three-manifold. We equip each copy of \(C\), indexed by \(w\in W\), with the metric \(h\) chosen in Section 2.4. Since \(h\) has the product form \(h_j+dt^2\) on each collar, these metrics fit together smoothly across the identified boundary components. In the signed collar coordinates, the resulting metric is again \(h_j+dt^2\). Since \(\Sigma\) is obtained by gluing finitely many copies of the compact manifold \(C\), it is compact, and the metric is therefore complete. Moreover, $\Sigma$ is connected, since crossing a
boundary component changes $w$ to $ws_j$ and the elements $s_j$
generate $W$. Left multiplication of the labels defines an isometric action of \(W\) on \(\Sigma\), given by \(w'\cdot[w,c]=[w'w,c]\) for \(w'\in W\).

Now define
\begin{equation}\label{eq:iota}
 \iota\colon\Sigma\longrightarrow M
 \ \ \text{ by } \ \ \iota([w,c])=q(w,c).
\end{equation}
This is well-defined because \(s_j\in K_j\), so \(q(w,c)=q(ws_j,c)\) whenever \(c\in F_j\). It is also equivariant
for the $W$-actions just described. Over the interior of each
chamber copy it is an isometric embedding. Recall that \(\Phi_j\) identifies the quotient over the \(j\)th collar with \(F_j\times E_j\). In these coordinates, \(\iota\) takes the form
\begin{equation}\label{eq:localimm}
 (\Phi_j\circ\iota\circ\eta_{w,j})(y,t)=(y,[w,t]_j).
\end{equation}
For $t<0$, this follows from
$[ws_j,-t]_j=[w,t]_j$, since $\chi_j(s_j)=-1$.
Thus the two half-collars map to the two halves of one real line
in the disk fiber. The map is smooth and has rank three at $t=0$
as well as away from zero. Its induced metric is $h_j+dt^2$, so
$\iota$ is an isometric immersion on all of $\Sigma$.

On the principal stratum the tangent image of $\iota$ is the
$C$-direction and is perpendicular to the orbit. At a point over
$F_j$, it is the sum of $T_yF_j$ and a real line in the disk fiber.
At the zero section both of these summands are perpendicular to
the tangent space of $G/K_j$. Since the identity copy of $C$
already meets every orbit, $\iota$ is a section. By the choice of the metric, each \(\{w\}\times\operatorname{Int}C\) is totally geodesic. The second fundamental form of \(\iota\) therefore vanishes on a dense open subset of \(\Sigma\) and hence everywhere. Thus \(\iota\) is a totally geodesic immersion.

There are no identifications in \(M\) over \(\operatorname{Int}C\), so \(\iota\) is injective on \(\iota^{-1}(M_{\mathrm{reg}})\). It therefore cannot factor through a covering \(\Sigma\to\Sigma'\) of degree greater than one, since such a factorization would give every point of \(\iota(\Sigma)\cap M_{\mathrm{reg}}\) at least two preimages under \(\iota\).

The behavior over $F_3$ is different. For $c\in F_3$, the points  $p=[e,c]$ and $p'=[r,c]$
of $\Sigma$ are distinct, since the gluing in
\eqref{eq:Sigma} identifies $e$ only with $s_3=r^2$ at this
boundary component. Their images agree, however, because
$r\in K_3$ and hence $q(e,c)=q(r,c)$. In the disk fiber over
$eK_3$, their local branches are
\[
 (y,t)\longmapsto(y,[e,t]_3) \ \ \text{ and } \ \
 (y,t)\longmapsto(y,[r,t]_3)=(y,[e,-it]_3).
\]
Here we have used $\chi_3(r)=-i$, which follows from
\eqref{eq:r}. At their common image point \((c,[e,0]_3)\) the tangent spaces are
therefore
\begin{equation}\label{eq:branches}
 T_cF_3\oplus\R
 \quad\text{and}\quad
 T_cF_3\oplus i\R
 \ \ \text{ in } \ T_cF_3\oplus\C.
\end{equation}
They intersect along \(T_cF_3\), and the two additional directions are orthogonal. Thus the two local branches meet with distinct tangent spaces.

The same argument determines all self-intersections.
For $c\in F_j$, a point $[w',c]$ maps to $q(w,c)$ precisely when
$w^{-1}w'\in W\cap K_j$. In \(\Sigma\), two such labels represent the same point exactly when they differ by right multiplication by \(s_j\), i.e. the gluing of
$\Sigma$ identifies exactly $w'$ and $w's_j$. Consequently
\begin{equation}\label{eq:multiplicity}
 \#\iota^{-1}(q(w,c))
 =\frac{|W\cap K_j|}{2}
 \ \ \ \text{ for } \ w\in W\ \text{ and }\  c\in F_j.
\end{equation}
By \eqref{eq:intersections}, this number is one over $F_1$ and
$F_2$ and two over $F_3$, since $r$ has order four. 
For each $c\in F_3$, the points $[w,c]\in\Sigma$ with $w\in W$
have exactly two distinct images under $\iota$, namely $q(e,c)$
and $q(s_1,c)$. Indeed, every element of $W$ has the form $r^k$
or $s_1r^k$, and right multiplication by $r^k\in K_3$ does not
change the image under $q$ over $F_3$. The entire self-intersection
set in $M$ is therefore
\[
 q(\{e\}\times F_3)\ \cup\ q(\{s_1\}\times F_3).
\]
Under the identification $S_3\cong F_3\times G/K_3$, these sets
are $F_3\times\{eK_3\}$ and $F_3\times\{s_1K_3\}$, so they
are disjoint embedded $2$-spheres. Every point on either sphere has
exactly two preimages under $\iota$. The tangent images of the
two branches are given by \eqref{eq:branches} on the first sphere
and by their translates under $s_1$ on the second.
This proves Theorem~\ref{thm:main}.

\begin{remark}
The self-intersections reflect a difference between the stabilizers
of a point of $\Sigma$ and its image in $M$. For $c\in F_3$, the
stabilizer of $[e,c]$ under the $W$-action is
$\langle s_3\rangle$, whereas the stabilizer of $q(e,c)$ under
the restricted $W$-action on $M$ is $\langle r\rangle$.
Thus the quarter-turn $r$ fixes the image point but exchanges
its two preimages in $\Sigma$.

In the Grassmann-bundle description of the geometric resolution
\cite{Lyt10}, these two branches are separated by their tangent
three-planes. In fact, the tangent lift
\[
 p\longmapsto d\iota_p(T_p\Sigma)\in\Gr_3(TM)
\]
is an embedding. It is an immersion because its composition
with the bundle projection is $\iota$, and it is injective
because distinct points with the same image under $\iota$
have distinct tangent images. Compactness of $\Sigma$ then
gives the embedding. The self-intersections occur when this
embedded lift is projected back to $M$.
\end{remark}

\end{document}